\pdfoutput=1
\documentclass{article}
\usepackage{graphicx} 

\usepackage{graphicx} 
\usepackage{amssymb}  
\usepackage{amsmath}
\newcommand{\ol}{\overline}

\newcommand{\ra}{\rightarrow}
\newcommand{\s}{\subseteq}
\newcommand{\es}{\emptyset}
\newcommand{\Ra}{\Rightarrow}
\newcommand{\LRa}{\Leftrightarrow}
\newcommand{\Q}{\mathbb{Q}}

\title{Enumerating all geodesics }

\author{Marcel Wild}

\begin{document}

\maketitle

\begin{quote}
A{\scriptsize BSTRACT}: {By {\it geodesic} we mean any sequence of vertices\\ $(v_1,v_2,...,v_k)$ of a graph $G$ that constitute a shortest path from $v_1$ to $v_k$. We propose a novel, 
natural  algorithm to enumerate all geodesics of $G$, and pit it (using Mathematica) against the standard procedure for the task.  The distance matrix $D(G)$ plays a crucial role in this. In fact, part of our article is devoted to survey its many uses in related tasks.}

\end{quote}

\section{Introduction}

\vspace{4mm}
Finding geodesics in graphs (which may also be directed and/or weighted) has a long history.
Up to one relevant exception, our article is a {\it survey} of these matters. Specifically it is structured along these lines:


\begin{itemize}
    \item[(OG)] What is the crispiest\footnote{In a minute we explain how "the crispiest" is to be understood.} way to find {\it one geodesic} between two fixed vertices? (Section 2)
    
    \item[(AG)] What is the crispiest way to find {\it all geodesics} between to fixed vertices?
    (Section 3)
    
    \item[(APAG)] What is the crispiest way to solve (AG) for {\it all pairs} of vertices?\\
    (Sections 4 an 6)
\end{itemize}

The  first instinct concerning (APAG) is to apply (AG) many times. However, since (AG) relies on depth-first search (DFS), whereas our method (called {\tt Natural-APAG}) avoids DFS, the "first instinct" looses out on {\tt Natural-APAG}. This is testified by numerical experiments in Sections 4 and 6.

\vspace{3mm}
Most everything in this article hinges on our assumption that the distance matrix $D(G)$ is known\footnote{Without these rose-colored glasses the (OG),(AG),(APAG)  landscape looks  rocky. Yet [S]  is an excellent guide to explore it.}. In fact, {\it given} $D(G)$, the author attempts to find the crispiest  way to solve each of (OG),(AG),(APAG). Roughly put "crispiest" means  that there are hardly faster methods when $D(G)$ is {\it known}. And if there are, they must be rather clumsy. 

All graphs are assumed to be  {\it connected}\footnote{This is just for convenience and no loss of generality. Alternatively we must always pick some connected component of our graph.}.
Furthermore, up to and including Section 4, our graphs will be undirected and unweighted.

Let us point out clearly the two  different aspects of our article. On the one hand it surveys, for non-experts, the benefits of the distance matrix $D(G)$. This mainly from a mathematical point of view, but with a glimpse at applications in the social and biological sciences. On the other hand, our article features an original, competitive method concerning one particular aspect of $D(G)$, namely the computation of all geodesics of $G$.

An effort was made to avoid excessive formalism, and rather use toy examples.
In fact, no (traditional) proof will be encountered. Quite a bit of additional information was put in footnotes in order not to interrupt the flow of thoughts.
Note that "iff" will be used as a shorthand for "if and only if".

\section{The distance matrix}

Although   the reader is assumed to be familiar with basic graph theory, we still need to fix some notation in Subsection 2.1. The distance matrix $D(G)$ is introduced in 2.2, and exploited in 2.3 to solve problem (OG).

\includegraphics[scale=1.1]{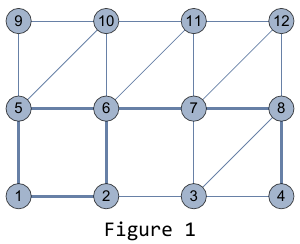}

\vspace{1cm}

{\bf 2.1} Consider the graph $G_1$ in Figure 1. It has vertex set $V_1=\{1,2,...,12\}$ and an edge-set of cardinality $|E_1|=24$. Formally, if $G=(V,E)$ then the elements of $E$, i.e. the {\it edges}, are unordered pairs of vertices. One says that vertices $x,y\in V$ are  {\it adjacent} if $\{x,y\}\in E$. A {\it walk} $W$ from a {\it source} $s\in V$ to a {\it sink} $t\in V$ is any sequence $(s,v_1,...,v_k,t)$ of vertices  such that $v_i,v_{i+1}$ are adjacent for all $0\le i\le k$ (upon setting $v_0:=s$ and $v_{k+1}:=t$). The number of involved edges, i.e. $\ell(W):=k+1$, is called the {\it length} of $W$. If we like to point out source and sink, we speak of $s-t$ walks. For instance $W_1:=(3,7,11,12,7,6)$ and $W_2:=(3,7,12,11,6)$ are both 3-6 walks. They have the same {\it underlying vertex-set} $V(W_1)=V(W_2)=\{3,6,7,11,12\}$, but $\ell(W_1)=5\neq 4=\ell(W_2)$.

A  walk without repeated vertices (a fortiori: without repeated edges) is called a {\it path}. Thus $W_2$  is a path but not $W_1$.
A {\it shortest} $s-t$ walk  is one of minimum length. Evidently it must be a path.
Such paths (for whatever $s,t$) will usually be called {\it geodesics} [H,p.324]. For instance $(3,7,6)$ is a geodesic. 
Notice that each subpath of a geodesic is itself a geodesic (why?).
Furthermore the following is easily verified. If $P_1$ and $P_2$ are geodesics with $V(P_1)=V(P_2)$ then\footnote{This fails when geodesics are defined in terms of {\it weight} instead of length, as happens in 5.2.} $P_1=P_2$.

The {\it distance}  between two vertices $s\neq t$ is the length $dist(s,t)$ of any $s-t$ geodesic. We further put $d(s,s):=0$.
Two questions arise: First, how to find {\it one}  $s-t$ geodesic? Second, how to get $dist(s,t)$? Settling the first establishes the second. Moreover, there seems to be no other way to settle the second. As seen in 2.3, this is false. 

\vspace{3mm}

{\bf 2.2} Fix some integer $k\ge 1$. Suppose that  $G$ has vertex set $V=\{1,2,...,n\}$ and that we know the number $b_{i,j}$ of walks from $i$ to $j$ which have length $k$.  Let $a_{i,j}$ be $1$ if the vertices $i,j$ are adjacent, and $0$ otherwise.
Suppose that $x,y$ are any fixed (not necessarily distinct) vertices. Then $b_{x,3}a_{3,y}$ is the number of walks of type $(x,...,3,y)$ where $(x,...,3)$ ranges over all (if any) length $k$ walks from $x$ to $3$. Either $b_{x,3}=0$ or $a_{3,y}=0$ can cause
$b_{x,3}a_{3,y}=0$. In any case, $b_{x,3}a_{3,y}$ counts the number of length $k+1$ walks from $x$ to $y$ whose second last vertex is $3$. By considering not just $z=3$ but all $z\in\{1,2,..,n\}$
it follows that the number $c_{x,y}$ of walks from $x$ to $y$ of length $k+1$ equals

$$(1)\quad c_{x,y}=b_{x,1}a_{1,y}+b_{x,2}a_{2,y}+\cdots+b_{x,n}a_{n,y}.$$

 The $n\times n$ matrix $A=A(G)$  with $(i,j)$-entry $A(i,j):=a_{i,j}$ is called the {\it adjacency matrix} of $G$. Similarly the  $n\times n$ matrices $B$ and $C$ are the ones with entries $B(i,j):=b_{i,j}$ and $C(i,j):=c_{i,j}$. It follows from (1) that $C=BA$. Furthermore, if $k=1$, then $B=A$. Thus for $k=1$ it holds that $C=A^2$. Hence for $k=2$ it holds that $B=A^2$, and so $C=BA=A^3$. By induction on $k$ one concludes that
 $A^{k+1}(x,y)$ is the number of $x-y$ walks of length $k+1$.
 
 This has pleasant consequences: Suppose $\alpha:=A^{k+1}(x,y)\neq 0$ but $A^k(x,y)=0$. Then necessarily $dist(x,y)=k+1$. Moreover, the set of all  $\alpha$  walks from $x$ to $y$ coincides with the set of all shortest $x-y$ paths, i.e. $x-y$ geodesics. Yet we have no idea how these  geodesics look like! 

 \vspace{3mm}

  The {\it distance matrix} $D=D(G)$ by definition has entries $D(i,j):=dist(i,j)$. By the above we can obtain $D$ by taking powers $A,A^2,A^3,..$ of the adjacency matrix until finally $A^\delta$ has zeros only on the diagonal. Using standard matrix multiplication this takes $O(n^3\cdot\delta)$ time. The upper bound $n$ for  $\delta$ is sharp. One may frown at the high complexity $O(n^3\cdot\delta)=O(n^4)$ to calculate  $D(G)$ in the "crisp way", but $O(n^4)$ only hurts for  large values of $n$. To assess what "large" may mean, calculating $D(G)$ the crisp way for some $G$ with 6000 vertices and 20000 edges took\footnote{More precisely, it took  7.6 seconds to calculate $A,A^2,..A^\delta=A^9$, but we didn't bother to replace the arising nonzero entries in the powers $A^k$ by the appropriate distance-values.} modest 7.6 seconds. 
  Admittedly,  for larger dimensions this method cannot finish in reasonable time; more on that in 6.3. 

  \vspace{3mm}
 
 {\bf 2.3} Suppose that $D$ is known. Given $s,t\in V$, some {\it single}  $s-t$ geodesic can then be found as follows. Say that by consulting $D$ we know that $dist(s,t)=k$. Among the {\it neighbors} (:=adjacent vertices) of $s$ there must be at least one $x_1$ with $d(x_1,t)=k-1$. (This kind of reduction of optimality to the optimality of substructures is the core of {\it dynamic programming}.)
  Using $D$ we can find $x_1$ efficiently.
  Likewise, among the  neighbors of $x_1$ we find some $x_2$ with $d(x_2,t)=k-2$, and so on until we have a path $(s,x_1,...,x_{k-1},t)$. This arguably is the crispiest\footnote{There is a way [S,p.88] to find a shortest $s-t$ path that avoids $D(G)$, but it is "overkill" in the sense that it necessarily produces a shortest $s-t$ path for {\it all} $t\in V$. The arising structure is called a {\it shortest-path-tree} rooted at $s$. Recall our survey assumes that $D(G)$ is given.} way
  to find a  $s-t$ geodesic. Therefore problem (OG) in Section 1 is settled. What about (AG)?

\section{Solving problem (AG)}

When one desires {\it all} geodesics  between fixed vertices $s\neq t$,
this could be reduced to (OG), i.e. by choosing $x_1,x_2,...$ in 2.3 in all possible ways. It isn't however obvious how to carry this out without regenerating previously obtained partial paths.

We will show how to organize the book keeping on the graph $G_1$ of Fig.1. 
To make  list (2) below more readible we relabelled for the time being $a:=10,\ b:=11,\ c:=12$.

\vspace{2mm}
{\bf 3.1} By inspection (in general: by consulting $D(G)$) one finds that $dist(1,c)=4$. In order to find {\it all} length 4 paths from $1$ to $c$, we will put partial paths in a stack, snapshots of which are shown in (2). 

Every path departing from $1$  either starts as $(1,2)$ or $(1,5)$. In shorthand: $12$ or $15$. This explains the first snapshot.
Momentarily leaving $15$ aside we turn to $12$ and build its extensions $126$ and $123$. They\footnote{We always add the extensions in "lexicographic" decreasing order. This is not strictly necessary, but this way in the end all shortest $s-t$ paths are output in lexicographic {\it increasing} order.} replace $12$ in the stack. In similar\footnote{Generally the path "Last In" (i.e. 123) is the path "First Out" that undergoes treatment (here: being replaced by 1238, 1237). Such a queue to handle any kind of "treatment" is called a LIFO-stack. It is well known that LIFO-stacks and depth-first-search (DFS) are two sides of the same coin.} fashion $(12,126,123)$ leads to $(15,126,1238,1237)$.
The extensions $12378$ and $1237b$ are duds since they do not extend to shortest paths (in general one needs $D(G)$ to weed out duds), but $1237c$ works. In fact, since  it happens to be a $1-c$ geodesic already, we render it in boldface. Generally geodesics are removed from the LIFO-stack and stored somewhere else. Continuing in this manner we get altogether seven (boldface) $1-c$ geodesics:

\begin{itemize}
    \item[(2)] $(15, 12)\  (15, 126, 123)\ (15,126,1238,1237)\ (15, 126, 1238, {\bf 1237c})$
    \item[] $ (15, 126, 1238)\ (15, 126, {\bf 1238c})\ (15, 126)\ (15, 126b, 1267)\ (15,126b, {\bf 1267c})$
    \item[] $(15, 126b)\ (15, {\bf 126bc})\ (15)\ (15a, 156)\ (15a, 156b, 1567)\ (15a,156b,{\bf 1567c})$
    \item[] $(15a,156b)\ (15a,  {\bf 156bc})\ (15a)\ (15ab)\ ({\bf 15abc})$
\end{itemize}

It is evident that this DFS method generalizes to the lexicographic enumeration of all $s-t$ paths of length at most $k$, for any fixed positive integer $k$. (When $k<dist(s,t)$ then finding nothing is still informative.) It is easy to see, and well known, that this type of DFS runs in time $O(N_{s,t}\cdot pol)$, where $N_{s,t}$ is the number of $s-t$ paths of length $\le k$, and $pol$ is some polynomial in $|V|$ and $|E|$ (which we omit to spell out).

\vspace{3mm}
{\bf 3.2} Some version of depth-first search likely lurks\footnote{Unfortunately, Mathematica's internal implementation is not publicly documented.} behind the (self-explanatory) Mathematica command ${\tt FindPath[s,t,k,All]}$. 
The output of ${\tt FindPath[1,12,5,All]}$ is shown in (3). 
Up to ordering and up to using $\{,\}$ instead of $(,)$, the first seven sets in (3) match the boldface entries above. The other 16 sets are the  $1-12$ paths of length 5. 

\vspace{3mm}

$(3)\quad \Big\{ \{1, 5, 10, 11, 12\}, \{1, 5, 6, 11, 12\}, \{1, 5, 6, 7, 12\}, \{1, 2, 6, 
  11, 12\},\\ \{1, 2, 6, 7, 12\},
  \{1, 2, 3, 8, 12\}, \{1, 2, 3, 7, 12\};\  
  \{1,  5, 10, 6, 11, 12\}, \{1, 5, 10, 6, 7, 12\},\\
  \{1, 5, 10, 11, 7, 12\},  \{1, 
  5, 9, 10, 11, 12\}, \{1, 5, 6, 11, 7, 12\}, \{1, 5, 6, 10, 11, 12\},\\
  \{1,  5, 6, 7, 11, 12\}, \{1, 5, 6, 7, 8, 12\}, \{1, 2, 6, 11, 7, 12\}, \{1, 2, 
  6, 10, 11, 12\},\\
  \{1, 2, 6, 7, 11, 12\},
  \{1, 2, 6, 7, 8, 12\}, \{1, 2, 3,
   8, 7, 12\}, \{1, 2, 3, 7, 11, 12\}, \\
   \{1, 2, 3, 7, 8, 12\}, \{1, 2, 3, 4, 
  8, 12\} \Big\}$

\section{Introducing {\tt Natural-APAG}}

 Recall from Section 1 that APAG is the problem to enumerate, for all  vertex-pairs\footnote{The attentive reader has noticed that it suffices to get all  $s-t$ geodesics only for $s<t$ because upon reversing them one gets all  $t-s$ geodesics. This will no longer apply in Section 6 where we deal with directed graphs.} $(s,t)$ all  $s-t$ geodesics. An obvious way to achieve this is to apply $\binom{n}{2}$ times (AG) from Section 3.
  We dare to challenge this approach.
 The technical details of our algorithm {\tt Natural-APAG} are postponed to Section 6 (where we also generalize the matter to weighted digraphs).

 \vspace{2mm}
 {\bf 4.1} Whatever they are,  {\tt Natural-APAG} has the inherent disadvantage of being set up in high-level Mathematica code, whereas {\tt FindPath} is a hard-wired Mathematica command. 
Nevertheless, as Table 1 indicates,  {\tt Natural-APAG}  prevails the more  the {\it fewer} and the {\it longer} the $(s,t)$-geodesics are on average. As to "fewer", it does not pay off to relaunch {\tt FindPath} $\binom{n}{2}$ times, only to collect each time a few geodesics,  while {\tt  Natural-APAG} runs uninterrupted throughout. As to "longer", {\tt FindPath} is based on depth first search and hence struggles more than {\tt Natural-APAG} to generate a thousand geodesics of cardinality 25 than  a thousand of cardinality 9.

This is testified by Table 1 where by definition  $GridGraph[n_1,n_2]$ is the direct product of a $n_1$-vertex path with a a $n_2$-vertex path. In plain words, $GridGraph[n_1,n_2]$ is visualized by a $n_1\times n_2$ grid-paper. Furthermore\\  $CircGraph[n,\{\alpha,\beta\}] $ has vertices $1,2,...,n$ and each  vertex $i$ is adjacent
(modulo) $n$ to the vertices $i-\alpha,\ i+\alpha,\ i-\beta,\ i+\beta$.

\vspace{3mm}

\begin{tabular}{|c|c|c|c|c|}
			
 Graph, $maximum\ distance$ & $\#$ geodesics& {\tt Nat-APAG} & iter. {\tt FindPath}    \\ \hline
             & &    &                                 \\ \hline
   $GridGraph[10,12]$, 20&  4'990'560 &97 sec &  18 sec    \\ \hline
   $GridGraph[5,24]$, 27&  1'470'400 &46 sec &  1545 sec    \\ \hline
   $CircGraph[100,\{1,2\}]$, 25&  35000 &0.8 sec &  16213 sec    \\ \hline
   $CircGraph[100,\{4,7\}]$, 9&  69000 &2.2 sec &  0.2 sec    \\ \hline
   
\end{tabular}

\vspace{2mm}

{\sl Table 1: Pitting high-level {\tt Natural-APASP} against hard-wired {\tt FindPath}}

\vspace{2mm}

 The distance matrix $D(G)$ is crucial for both algorithms and was computed (whatever way) in lightning speed by  the Mathematica command\\
 {\tt GraphDistanceMatrix[G]}; more details in 6.3.

 
 \vspace{3mm}
{\bf 4.2} As observed in 2.1, each subpath of a geodesic $P$ is itself a geodesic.  In other words, if $X\s V$ is the set of vertices underlying $P$ then

\begin{itemize}
\item[(4)]{\it For all $a,b\in X$ there is a geodesic that lies in $X$ and contains $a,b$.}
\end{itemize}

Generally, if $X\s V$ is {\it any} subset satisfying (4) then $X$ is called a {\it metric} subset. It turns out [W] that, once all geodesics are known, all metric subsets can be enumerated in compressed fashion (using don't-care symbols). If $X\s V$ is metric then the induced graph $G[X]$ is not just (ordinarily) connected but "supremely" connected. It is shown in [W] that, like the connected sets, also the metric sets can  be enumerated in output polynomial time. Numerical experiments are provided as well.

Why bother about metric subsets? Enumerating all (ordinarily) connected induced subgraphs has received considerable attention recently. The key word here is {\it Community Detection}\footnote{We recommend [GN] as an introduction  that features many concrete examples. Use GoogleScholar to browse the more than 20'000 articles citing [GN] and find the ones that  specifically focus on the enumeration of all connected sets. See also [P].}  in social or biological networks. Directing attention from connected to metric subsets  delivers tighter communities and often cuts the  number of enumerated entities from the millions to the thousands.


 \section{Digraphs, weighted or not}

 In 5.1 we fix notation and note that for strongly connected {\it digraphs} $G$ our crisp way to obtain $D(G)$ persists. In 5.2 we allow our digraphs to be weighted. Unfortunately this destroys the crisp way towards $D(G)$. An alternative way is sketched in 5.3.
 
 \vspace{2mm}
 {\bf 5.1} Let us turn to directed graphs, aka {\it digraphs}. Hence by definition each edge (now called {\it arc}) has a direction\footnote{Undirected graphs can be viewed as those digraphs where each arc $(i,j)$ has the companion $(j,i)$.}. Formally the arc-set of a digraph $G=(V,Arc)$ is a set $Arc\s V\times V$ of {\it ordered} vertex-pairs $(i,j)$. For technical reasons we forbid loops, i.e.  arcs of type $(i,i)$.
 A {\it directed $s-t$ walk}  is any sequence $(s,v_1,...,v_k,t)$ of vertices $s:=v_0,v_1,...,v_k,v_{k+1}:=t$ such that $(v_i,v_{i+1})\in Arc$  for all $0\le i\le k$. A {\it directed path (dipath)} $P$ is a directed walk that doesn't repeat vertices. The {\it length} $\ell(P)$ is the number of arcs in $P$. By definition $dist(s,t)$ is the length of a shortest  $s-t$ dipath $P$. Alternatively we speak again of $s-t$ geodesics. Furthermore $dist(s,s):=0$ and $dist(s,t):=\infty$ if $s\neq t$ and there is no  $s-t$ dipath. The latter never happens
 iff $G$ is {\it strongly connected} in the usual sense [S,p.32] that there is some $s-t$ dipath for all $s,t\in V$. 

\vspace{2mm}
{\bf 5.1.1}
 We again denote by $A=A(G)$ the matrix whose $(i,j)$-entry $a_{i,j}$ is $1$ if $(i,j)\in Arc$, and $0$ otherwise. Similarly $D(G)$ is the matrix whose $(i,j)$-entry is $dist(i,j)$. As opposed to Section 2, usually $A$ and $D(G)$ are not symmetric. We henceforth assume that $G$ is strongly connected in order to have all matrix entries $dist(i,j)\neq\infty$.   
  It is easy to see that $D(G)$ can, as in Section 2, be obtained by calculating $A,A^2, A^3,...$. Similarly, our solution in Section 3 for (AG) carries over from connected graphs  to strongly connected digraphs.

\vspace{2mm}
{\bf 5.2}  Let us turn to weighted digraphs. We thus have a function $w:Arc\to \Q$ which to each arc $(i,j)$ assigns its {\it weight}\footnote{Ordinary (unweighted) digraphs can be viewed as weighted digraphs all of whose arcs have weight 1. In this case $w=\ell$, i.e. weight and length are the same thing.} $w(i,j)$. The weight $w(P)$ of a dipath $P$ is defined as the sum of the weights of its arcs. Do not confuse $w(P)$ and $\ell(P)$.
We define $dist_w(s,t)$ as the weight of a $s-t$ {\it geodesic}, this being a 
$s-t$ dipath of minimum weight.  
By definition the $(i,j)$-entry of the {\it distance-matrix} $D_w(G)$ is $dist_w(i,j)$. 
Easy examples show that   shortest $s-t$ dipaths 
need not be  $s-t$ geodesics. 

Therefore the following is no surprise. 
 Let $A_w(G)$ be the adaption of the adjacency matrix $A(G)$ of 5.1.1 where each entry $\ell(i,j)=1$  is replaced by $w(i,j)$. Then, unfortunately,   $D_w(G)$ {\it can no longer} be obtained by evaluating powers of $A_w(G)$.

\vspace{3mm}
{\bf 5.3} 
Here comes a smart (and crisp) way to calculate $dist_w$, and hence $D_w(G)$.
Interestingly, it even
works  in the presence of negative weights, provided that $G$ has {\it no negative directed circuits}.
Closely following [S,p.110] fix an arbitrary ordering $v_1,...,v_n$ of the vertex-set $V$ of $G$. For $s,t\in V$ and $0\le k\le n$ define

\begin{itemize}
    \item[(5)] $d_k(s,t):=$ minimum weight of an $s-t$ directed walk using only vertices in $\{s,t,v_1,...,v_k\}$.
\end{itemize}

\noindent
Then it holds that $d_0(s,t)=w(s,t)$ if $(s,t)\in Arc$, while $d_0(s,t)=\infty$ otherwise. Moreover

$$(6)\quad d_{k+1}(s,t)=min\big\{d_k(s,t),\ d_k(s,v_{k+1})+d_k(v_{k+1},t)\big\}$$

\noindent
for all $s,t\in V$ and $k<n$. It is easy to see that $dist_w=d_n$ and that $d_n$ can be calculated in time $O(n^3)$. 


\section{{\tt Natural-APAG} da capo}

In this Section weights must be positive {\it and}  integer, by reasons that will be evident soon. Thus our digraph $G$ comes with a weight function $w:Arc\to\{1,2,\ldots\}$. Using a lush toy example (6.1) we are going to disclose the inner workings of {\tt Natural-APAG}, which was numerically evaluated for ordinary graphs already in Section 4. {\tt Natural-APAG} actually comes in two variants: backwards (introduced in 6.1) and forwards. The complexities of both variants are determined in 6.2 and 6.3 respectively. Section 6.4 features comprehensive numerical experiments that go beyond the ones (for ordinary graphs) in Table 1.
In 6.5 we generalize geodesics to chordless paths and hint at real-life applications.

\vspace{3mm}
{\bf 6.1}
Let us illustrate {\tt Natural-APAG} on the weighted digraph $G_2$ in Figure 2. We write $Geo[k]$ for the set of  geodesics  of  weight ({\it not} length)  $k$. The overall strategy is to recursively calculate $Geo[k]$ based on $Geo[1],...,Geo[k-1]$, until  all geodesics have been found.

\vspace{3mm}

\includegraphics[scale=1.15]{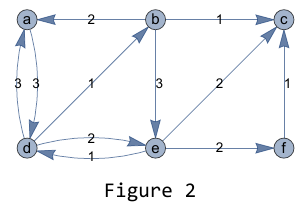}

\noindent
We thus start with all weight 1 arcs, i.e. with

$$Geo[1]=\{bc,db,{\bf ed},fc\},$$

\noindent
where e.g. $db$ is a shorthand for the dipath $(d,b)$. (The entry {\bf ed} is boldface to aid the discussion of $Geo[4]$ below.) Generally each geodesic  $P\in Geo[k]$, which isn't a single arc, is a unique
extension\footnote{The converse fails.  For instance, although $db\in Geo[1]$ and $w(be)=3$, it is false that $dbe\in Geo[4]$ since $dist_w(d,e)=2$. Therefore the graph-distance-matrix $D_w(G)$ is mandatory, but for $G_2$ inspection will do.} on the  right of some geodesic $Q$ with $w(Q)<w(P)$. Hence $Q\in Geo[h]$ for some $h<k$. Specifically, $k-\mu\le h<k$, where $\mu$ is the maximum occurring  arc-weight. In our case $\mu=3$. 
 
When we build $Geo[k]$ by scanning some of its forerunners $Geo[h]$ we will proceed from $h=k-1$ to  $h=k-2$ to $  h=k-3$, as long as these integers are $\ge 0$. The occurrence of $h=0$ amounts to the "extension" of a vertex to an arc. Let's go:

\begin{itemize}
    \item[] $Geo[2]=\{dbc, {\bf edb;}\ ba,{\bf de},ec,ef\}$
    \item[] $Geo[3]=\{--{\bf ;}\ dba{\bf ;}\ {\bf ad, be}, da\}$
    \item[] $Geo[4]=\{adb,bed{\bf ;}\ edba,def{\bf ;}\ eda\}$
    \item[] $Geo[5]=\{adbc{\bf ;}\ ade, bef\}$
    \item[] $Geo[6]=\es$
\end{itemize}

\noindent
Consider say $Geo[4]$. By scanning $Geo[3]$ one-by-one we find that the head $a$ of $dba$ has no outgoing arc of weight 1. But ${\bf ad}\in Geo[3]$ has $w(db)=1$.
 Additionally it happens that $dist_w(a,b)=4$. Hence $ad$  yields the first member $adb$ of $Geo[4]$. Also ${\bf be}\in Geo[3]$ {\it works} in the twofold sense that its head {\it has} outgoing arcs of weight 1, {\it and} for at least one of these "sensors" the  distance matrix gives "thumbs up". In contrast,  $da\in Geo[3]$ does not work (since $da$ does not even have sensors). Next we scan $Geo[2]$ and find that ${\bf edb}$ and ${\bf de}$ work. As to $de$, it has the two sensors $ec$ and $ef$
 (i.e. both have weight 2), but only the latter yields a geodesic $def$ (why?). Finally, in $Geo[1]$ only ${\bf ed}$ works.
 To summarize, exactly the boldface members in $Geo[1]\cup Geo[2]\cup Geo[3]$ contribute to $Geo[4]$.

We leave it to the reader to verify the claimed content of $Geo[5]$ and $Geo[6]$.
Notice that the emptiness of $Geo[6]$ would imply the emptiness of $Geo[7]$ in the unweighted case. Not here: $Geo[7]=\{adef\}$, but indeed
$Geo[8]=Geo[9]=Geo[10]$.
Generally, it follows from $Geo[k-\mu]=Geo[k-\mu+1]=\cdots =Geo[k-1]=\es$ that $Geo[h]=\es$ for all $h\ge k-\mu$. Hence in our case $Geo[h]=\es$ for all $h\ge 8$.

\vspace{3mm}
{\bf 6.1.1}
For any digraph with weight function $w:V\to\{1,2,...\}$ and vertices $s\neq t$ let $Geo[s,t]$ be the set of all $s-t$ geodesics. Then  $Geo[s,t]$ (possibly empty\footnote{This happens iff the digraph is not strongly connected.}) must be a subset of some $Geo[k]$. Put another way, each $Geo[k]$ can be written uniquely (up to order) as

$$(7)\quad Geo[k]=Geo[s_1,t_1]\uplus Geo[s_2,t_2]\uplus\cdots\uplus Geo[s_\ell,t_\ell],$$

\noindent
where we take the case $\ell=0$ to mean $Geo[k]=\es$. In our example one checks that  $Geo[c,x]=\es$ for all $x\neq c$, 
$Geo[f,x]=\es$ for all $x\not\in\{f,c\}$, $Geo[d,a]=\{da,dba\},\ Geo[e,a]=\{eda,edba\}$, and all other $Geo[s,t]$ are singletons.

\vspace{3mm}
{\bf 6.2} As to the complexity of {\tt Natural-APAG}, let $pol$ be some\footnote{Recall that $D(G)$ is given,  but it depends on the particular data structure on how fast one can search within $D(G)$; these details are not relevant here. }  low degree polynomial in the variables
$n=|V|$ and $m=|Arc|$ that bounds the time for doing the following to each given geodesic $P=(a,...,v)$ and integers $\alpha,\beta\ge 1$: For all $(v,b)\in Arc$ with $w(v,b)=\alpha$ check whether $dist_w(a,b)=\beta$.

Finding any one of the $N$ many geodesics $Q$ of $G$ requires at most the following amount of time. If $Q\in Geo[k]$ then for each geodesic $P$ in
$$Geo(Q):=Geo[k-\mu]\cup Geo[k-\mu+1]\cup\cdots\cup Geo[k-1]$$
decide whether $P$ "works" (in the sense of the toy example). Hence it takes $\le|Geo(Q)|\cdot pol\le Npol$ time to find $Q$. Therefore {\it all} $N$ geodesics of $G$ can be calculated in $O(N^2pol)$ time.


\vspace{3mm}
{\bf 6.3} Let us call the previously discussed algorithm the {\it backward-variant} of  {\tt Natural-APAG}. The {\it forwards-variant} is defined as follows. For starters, the set $Geo[1]$ is defined as in 6.1 and additionally {\it each} arc $(v,w)$ of $G$ is put in the appropriate box; i.e. if $(v,w)$ has weight (not length) $k$, it is put into $Geo[k]$.

By induction (anchored in $k=1$) assume that for some fixed $k\ge 1$  all boxes $Geo[j]\ (1\le j\le k)$
have been finalized. We process $Geo[k]$ as follows. For each geodesic $(a,...,b)\in Geo[k]$ scan {\it all} its sensors $(b,c)$. For each $(b,c)$ ask the distance-matrix whether $(a,...,b,c)$ is a geodesic. If yes, let $\ell$ be the weight of $(b,c)$ and put\footnote{The box $Geo[k+\ell]$ may or may not have been empty beforehand. } $(a,...,b,c)$ into $Geo[k+\ell]$. Having handled $Geo[k]$ as described, it is clear that all boxes $Geo[j]\ (1\le j\le k+1)$ are now finalized (and perhaps some more).
It is clear that the forwards variant of {\tt Natural APAG} also has complexity $O(N^2pol)$.

Let $N_{s,t}$ be the number of geodesics of $G$ whith startpoint $s$ and endpoint $t$. Using $D(G)$ it is trivial to decide whether
$N_{s,t}=0$. If not then by 3.1 (which extends to the weighted case) it takes time $O(N_{s,t}\cdot pol)$ to enumerate these particular geodesics. Because $N$ is the sum of all numbers $N_{s,t}\quad (1\le s,t\le n,\  s\neq t)$ it follows that the standard
procedure for APAG enumerates all $N$ geodesics of $G$ in time $O(N\cdot pol')$. Thus "on paper" the standard way beats {\tt Natural APAG}. However in practise (see 6.4) the picture is more subtle.

\vspace{3mm}
{\bf 6.4} In many instances of Table 2 the forwards-variant (fv) beats the backwards-variant (bv). In detail, Table 2 is set up as follows. For chosen numbers $n,m,maxwgt$ we create at random a directed and weighted digraph $G$ with $n$ vertices, $m$ arcs, and maximum\\ arcweight $maxwgt$.
Then the distance-matrix $D(G)$ is computed\footnote{Its computation time does not enter the CPU-time of any algorithm; even for the digraph $G$ with 12000 vertices computing the 144'000'000 entries of $D(G)$ only took 33.5 seconds.}, the maximum occuring distance $maxdis$ is recorded, and three algorithms are unleashed upon $G$. While all of them agreed upon the number of geodesics, the respective CPU times (in seconds)  often differed a lot.

In the first five instances $n=300,m=2000$ are fixed. As $maxwgt$ increases from $1$ to $1000$, the algorithm that suffers most is the backwards-version of {\tt Natural-APAG}. For example 229 versus 2 seconds for bv versus fv in the last instance.

In the next four instances $m=15000$ and $maxwgt=20$ are fixed. As $n$ increases from $1000$ to $5000$, iterating {\tt  FindPath}  falls behind. This is why. In the
$(5000, 15000,20)$-instance {\tt FindPath} must be applied to $5000\cdot4999=24'995'000$ vertex-pairs $(v_1,v_2)$ but on average there are only few\footnote{We did not check whether our random digraphs were strongly connected. If not, and there is no dipath from $v_1$ to $v_2$, this is  fast detected since $dist_w(v_1,v_2)=\infty$. If the (1000,15000)-digraph was strongly connected then on average only $23684428/24995000=0.95$  $(v_1,v_2)$-geodesics exist.} geodesics to be found
for each pair $(v_1,v_2)$. Furthermore these geodesics can have a weight up to 183. Already in Section 4 we observed that long geodesics are more painful for {\tt FindPath} than for {\tt Natural-APAG}. Likewise observe that as $maxwgt$ increases from 33 to 183, forwards {\tt  Natural-APAG} compares more and more favorably to the backwards variant. 

Similarly for the three digraphs with $n=6000,8000,12000$ and $maxdis=222,153,232$  both  backwards  {\tt  Natural-APAG}, and even more so {\tt FindPath}, would have  severly trailed  forwards  {\tt  Natural-APAG}. Interestingly, upon erasing weights, i.e. replacing (6000,70000,100) by (6000,70000,1), the number of geodesics explodes but the difference between bv and fv decreases.

The last four instances are of type $(1000, m,20)$. As the number of edges $m$ increases from $15000$ to $200'000$ (i.e. $G$ goes from sparse to dense), the performance pattern is not so clear-cut. It seems that when $m$ further increases the winner will be the {\it backwards} version of {\tt  Natural-APAG} but the runner-up is debatable. If 1000 and 20 become  other fixed values $n$ and $maxwgt$, and $m$ increases, it may well be that the pattern remains.

\vspace{5mm}

\begin{tabular}{|c|c|c||c|c|}
			
  $(n,m,maxwgt)\Ra maxdis$ & $\#$ geodesics & {\tt Nat-APAG}  & iter. {\tt FindPath}  
  \\ \hline
             & &    &                                 \\ \hline
   $(300, 2000,1)\Ra 6$&  261'518 &4.4/5.9 s &  1.8 s    \\ \hline
    $(300, 2000,4)\Ra 16$&  146'304 &3.6/3.3  s & 2.5 s   \\ \hline
    $(300, 2000,20)\Ra 60$&  101'019 &6.1/2.0 s  &  5.2 s  \\ \hline
     $(300, 2000,100)\Ra 275$&  92841 &22/2.1 s  &  3.4 s  \\ \hline 
      $(300, 2000,1000)\Ra 2763$&  89882 &229/2.0 s  &  3.7 s \\ \hline 
      
      
   & &    &                                 \\ \hline
  $(1000, 15000,20)\Ra 33 $&  1'340'556  & 100/61 s  &  137 s  \\ \hline
 $(2000, 15000,20)\Ra 68 $&  4'730'208  &760/297 s  &  4373 s  \\ \hline
    $(3600, 15000,20)\Ra 133$&  13'661'827 & 995/378 s  &  13436 s  \\ \hline
    $(5000, 15000,20)\Ra 183$&  23'684'428 & 3730/1124 s  &  31284 s  \\ \hline
 & &    &                                 \\ \hline
 $(6000, 70000,1)\Ra 6 $&  147'259'457  & 11663/9416 s  &  67984 s   \\ \hline
  $(6000, 70000,100)\Ra 222 $&  38'031'312  & ---/1555 s  &  ---   \\ \hline
   $(8000, 150'000,100)\Ra 153 $&  69'822'435  & ---/4535 s  &  ---   \\ \hline
     $(12000, 150'000,100)\Ra 232 $&  152'740'860  & ---/16177 s  & ---  \\ \hline
      & &    &                                 \\ \hline
  $(1000, 15000,20)\Ra 33 $&  1'340'556  & 100/61 s  &  137 s  \\ \hline
    $(1000, 40000,20)\Ra 14 $&  1'809'947  & 84/198 s  &  398 s  \\ \hline
      $(1000, 100'000,20)\Ra 8 $&  2'761'368  & 419/2471 s  &  839 s  \\ \hline
        $(1000, 200'000,20)\Ra 5 $&  4'350'787  & 153/4346 s  &  3268 s  \\ \hline   
\end{tabular}

\vspace{3mm}
{\bf 6.4.1} Here comes a more technical but crucial detail. Certain boxes $Geo[k]$ can end up with  millions of geodesics. Since adding new geodesics {\tt geod} to a  large crowd of geodesics (using the Mathematica command {\tt AppendTo[Geo[k]],geod]}) costs time, it pays to subdivide each box into $\alpha$ sub-boxes; not too few but also not too many. It turned out that $\alpha$ should be chosen  such  that each sub-box contains between 100 and 1000 geodesics. This entails that even $\alpha=100'000$ can occur (though $\alpha=1'000'000$ would be counter-productive).



\vspace{5mm}
\section{References}

\begin{itemize}

 \item[ {\bf GN} ] M. Girvan, M.E.J. Newman, Community structure in social and biological networks,
   Proceedings of the national academy of sciences 99.12 (2002): 7821-7826.

\item[ {\bf P} ] S. Pettie, All pairs shortest paths in sparse graphs, Encyclopedia of Algorithms, Springer 2008.

\item[ {\bf S} ]  A. Schrijver, Combinatorial Optimization, Springer Verlag 2003.

 \item[ {\bf W} ]  M. Wild, Compression with wildcards: All induced metric subgraphs, arXiv:2409.08363.

\end{itemize}

\end{document}